\documentclass[11pt]{article}

\usepackage{amsmath,graphicx}
\usepackage{epsfig}
\usepackage{graphicx}
\usepackage{psfrag}
\usepackage{amsmath}
\usepackage{amssymb,subfigure}
\usepackage{amsthm}
\usepackage{mathrsfs}
\usepackage{color}
\usepackage{pstricks,pst-node,pst-text,pst-3d}
\usepackage{pst-plot}

\def\R{\mathbb R}

\def\N{\mathbb N}
\def\tmu{\tilde{\mu}}

\def\M{\mathcal{M}}

\def\ub{\tilde{u}}

\def\epsilon{\varepsilon}
\def\ds{\displaystyle}

\def\ph{\phantom{1}}

\numberwithin{equation}{section}

\newcommand{\be}{\begin{equation}}
\newcommand{\ee}{\end{equation}}
\newcommand{\bes}{\begin{equation*}}
\newcommand{\ees}{\end{equation*}}
\newcommand{\baa}{\begin{array}}
\newcommand{\eaa}{\end{array}}
\newcommand{\ba}{\begin{eqnarray}}
\newcommand{\ea}{\end{eqnarray}}
\newcommand{\carre}{\hfill$\Box$\par\addvspace{4mm}}

\newtheorem{lemma}{Lemma}[section]
\newtheorem{theorem}[lemma]{Theorem}

\title{Uniqueness from pointwise observations in a multi-parameter inverse problem}

\author{Michel Cristofol$^{\hbox{\small{a}}}$, Jimmy Garnier$^{\hbox{\small{a,b}}}$, Fran\c cois Hamel$^{\hbox{\small{a,c,d}}}$ and Lionel Roques$^{\hbox{\small{b}}}$  \\
\footnotesize{$^{\hbox{a }}$Aix-Marseille Universit\'e, LATP, Facult\'e des Sciences et Techniques}\\
\footnotesize{Avenue Escadrille Normandie-Niemen, F-13397 Marseille Cedex 20, France}\\
\footnotesize{$^{\hbox{b }}$UR 546 Biostatistique et Processus Spatiaux, INRA, F-84000 Avignon, France}\\
\footnotesize{$^{\hbox{c }}$Helmholtz Zentrum M\"unchen, Institut f\"ur Biomathematik und Biometrie}\\
\footnotesize{Ingolst\"adter Landstrasse 1, D-85764 Neuherberg, Germany}\\
\footnotesize{$^{\hbox{d}}$ Institut Universitaire de France}}

\begin{document}

\maketitle

\begin{abstract}
In this paper, we prove a uniqueness result in the inverse problem of determining several non-constant coefficients of one-dimensional reaction-diffusion equations. Such reaction-diffusion equations include the classical model of Kolmogorov, Petrovsky and Piskunov as well as more sophisticated models from biology. When the reaction term contains an unknown polynomial part of degree $N,$ with non-constant coefficients $\mu_k(x),$ our result gives a sufficient condition for the uniqueness of the determination of this polynomial part.
This sufficient condition only involves pointwise measurements of the solution $u$ of the reaction-diffusion equation and of its spatial derivative $\partial u / \partial x$ at a single point $x_0,$ during a time interval $(0,\varepsilon).$ In addition to this uniqueness result, we give several counter-examples to uniqueness, which  emphasize the optimality of our assumptions.
Finally, in the particular cases $N=2$ and $N=3,$ we show that such pointwise measurements can allow an efficient numerical determination of the unknown polynomial reaction term.
\end{abstract}

\noindent{\it Keywords\/}: reaction-diffusion $\cdot$  inverse problem $\cdot$ multi-parameter $\cdot$ heterogeneous media $\cdot$ uniqueness

%
\section{Introduction}
%

Reaction-diffusion equations arise as models in many fields of mathematical biology \cite{Mur02}. From morphogenesis \cite{Tur52} to population genetics \cite{Fis37,kpp} and spatial ecology
 \cite{Shikaw97,Ske51,Tur98}, these partial differential equations benefit from a well-developed mathematical theory.

In the context of spatial ecology, single-species reaction-diffusion models generally deal with polynomial reaction terms. In a one-dimensional case, and if the environment is
supposed to be homogeneous they take the form:
\begin{equation}\label{eq_Intro1}
\frac{\partial u}{\partial t}-D  \frac{\partial^2 u}{\partial x^2}=P(u),
\end{equation}
where  $u=u(t,x)$ is the population density at time $t$ and space position $x$ and $D>0$ is the diffusion coefficient. The function $P,$ which stands for the growth of the population,
is a polynomial of order $N\ge 1$. The Fisher-Kolmogorov, Petrovsky, Piskunov (F-KPP) equation is the archetype of such models. In this model, we have $P(u)=\mu \, u - \gamma u^2$.
The constant parameters $\mu$ and $\gamma$ respectively correspond to the intrinsic growth rate and intraspecific competition coefficients. In this model, the lower the population density
$u$, the higher the per capita growth rate $P(u)/u.$ More complex models can involve polynomial nonlinearities of higher order. Examples are those taking account of an Allee effect.
This  effect occurs when the per capita growth rate $P(u)/u$ reaches its maximum at a strictly positive population density and is known in many species \cite{All38,Den89,VeiLew96}.
A typical example of reaction term involving an Allee effect is \cite{KeiLew01,LewKar93,RoqRoq08}:$$P(u)=r u\left(1-u\right) \left(u-\rho\right),$$with $r>0$
and $\rho \in (0,1)$. The parameter $\rho$ corresponds in that case to the ``Allee threshold" below which the growth rate becomes negative.

In the previous examples, the reaction terms were assumed to be independent of the space variable. However, real world is far from begin homogeneous. In order to take the heterogeneities into account, models have been adapted and constant coefficients have been replaced by space or time
dependant functions. In his pioneering work, Skellam \cite{Ske51} (and later, Shigesada, Kawasaki and Teramoto \cite{ShiKaw86}) mentioned the following extension of the F-KPP model to heterogeneous environments:
\begin{equation}\label{eq_IntroSke}
\frac{\partial u}{\partial t}-D  \frac{\partial^2 u}{\partial x^2}=\mu(x) \, u - \gamma(x) u^2.
\end{equation}
Here, the values of $\mu(x)$ and $\gamma(x)$ depend on the position $x$. For instance, regions of the space associated with high values of $\mu(x)$ correspond to favorable regions,
whereas those associated with low or negative values of $\mu(x)$ correspond to unfavorable regions. As emphasized by recent works, the precise spatial arrangement of these regions plays
a crucial role in this model, since it controls  persistence and spreading of the population \cite{BerHamRoq05a,CanCos03,ElsHamRoq09,RoqChe07,RoqHam07,RoqSto07,Shikaw97}. Models
involving an Allee effect can be extended as well to heterogeneous environments, as in \cite{HamFay10,RoqRoq08},
where the effects of spatial heterogeneities are discussed for models of the type:
\begin{equation}
\frac{\partial u}{\partial t}=D\frac{\partial^2 u}{\partial x^2}+r(x) \, u \left[(1-u) (u-\rho(x))+\nu(x)\right].
\label{bistable_heter}
\end{equation}
We also refer to \cite{MatNak06} for an analysis of propagation phenomena related to a reaction-diffusion model with an Allee effect in infinite cylinders having undulating boundaries.

In this paper, we focus on reaction-diffusion models with more general heterogeneous nonlinearities:
\begin{equation}
\frac{\partial u}{\partial t}=D\frac{\partial^2 u}{\partial x^2}+\sum_{k=1}^N\, \mu_k(x) \, u^k+g(x,u), \hbox{ for }t>0, \ x\in(a,b),
\label{eq_ini}
\end{equation}
for some interval $(a,b)$ in $\R.$

Since the behavior of such models depends on the precise form of the coefficients,  their empirical use  requires an accurate knowledge of the
coefficients. Unfortunately, in applications, the coefficients cannot be directly measured since they generally result from intertwined  effects of several factors.
Thus, the coefficients are generally measured through the density $u(t,x)$ \cite{SouHel08}. From a theoretical viewpoint, if $u(t,x)$ is measured at any time $t\ge 0$ and at
all points $x$ in the considered domain, all the coefficients in the model can generally be determined.  However, in most cases, $u(t,x)$ can only be measured in some -- possibly small --
subregions of the domain $(a,b)$ \cite{Wik03}.  For reaction-diffusion models as well as for many other types of models, the determination of the coefficients in the whole domain
$(a,b)$ bears on  inference methods which consist in comparing the solution of the model with hypothetical values of coefficients
$\tilde \mu_k,$ with the measurements on the subregions \cite{SouNeu09}. The underlying assumption behind this inference process is that there is a one-to-one and onto relationship between the value of the
solutions of the model over the subregions and the space of coefficients. This assumption is of course not true in general.

In this paper, we obtain uniqueness results for the coefficients $\mu_k(x)$, $k=1,\ldots, N$, based on localized measurements of the solution $u(t,x)$ of (\ref{eq_ini}).
The major differences with previous works dealing with comparable uniqueness results are (1) the size of the regions where $u(t,x)$ has to be known in order to prove uniqueness,
(2) the number of parameters we are able to determine, and~(3) the general type of nonlinearity we deal with.

Uniqueness of the parameters, given some values of the solution, corresponds to an inverse coefficient problem, which is generally dealt with -- for such reaction-diffusion equations --
using the method of Carleman estimates \cite{BukKli81,KliTim04}. This method provides Lipschitz stability, in addition to the uniqueness of the coefficients. However, this method requires,
 among other measurements, the knowledge of the density $u(\theta,x)$ at some time $\theta$ and for all $x$ in the domain $(a,b)$ (see \cite{BelYam06,BenCri08,CriGai06,ImmYam98,YamZou01}).   The uniqueness of the couple $(u,\mu(x))$  satisfying
the equation (\ref{eq_IntroSke}) given such measurements has been investigated in a previous work \cite{CriRoq08}, in any space dimension.

In a recent work, Roques and Cristofol \cite{RoqCri10} have proved the uniqueness of the coefficient $\mu(x)$ in \eqref{eq_IntroSke} when $\gamma(x)$ is known under the weaker assumption
that  the density $u(t,x_0)$ and its spatial derivative $\ds{\frac{\partial u}{\partial x}}(t,x_0)$ are known at a point $x_0$ in $(a,b)$ for all $t\in (0,\varepsilon)$
and that the initial density~$u(0,x)$ is known over $(a,b).$ This result shows that the coefficient $\mu(x)$ is uniquely determined in the whole domain $(a,b)$ by the
value of the solution $u(t,x)$ and of its spatial derivative at a single point~$x_0.$ The present work extends this result to the case of several coefficients $\mu_k(x)$, $k=1,\ldots, N:$
given any point $x_0$ in $(a,b)$ we establish a uniqueness result for the $N-$uple $(\mu_1(x),\ldots,\mu_N(x))$ given measurements of the $N$ solutions $u(t,x)$ of (\ref{eq_ini})
and of their first spatial derivatives in $(0,\varepsilon)\times \{x_0\}$, starting with $N$ nonintersecting initial conditions.

\section{Hypotheses and main result \label{sec:res}}

Let $(a,b)$ be a bounded interval in $\R$. We consider, for some $T>0$, the problem:
\begin{equation}\left\{\baa{l}\label{eq:multi_coef}
\ds\frac{\partial u}{\partial t}-D  \frac{\partial^2 u}{\partial x^2}=\sum_{k=1}^{N}{\mu_k(x)u^k}+g(x,u),\ \ t\in (0,T),\ \ x\in(a,b),\vspace{3pt}\\ [0.2 cm]
\ds\alpha_1 u (t,a) - \beta_1 \frac{\partial u}{\partial x} (t,a)=0,\ \ t>0,\vspace{3pt} \\ [0.2 cm]
\ds\alpha_2 u (t,b) + \beta_2 \frac{\partial u}{\partial x} (t,b)=0,\ \ t>0,\vspace{3pt}\\ [0.2 cm]
\ds u(0,x)=u^0(x),\ \ x\in(a,b),
\eaa
\right. \tag{$\mathcal{P}^{u^0}_{(\mu_k)}$}
\end{equation}
for some $N\in\N^*,$ and for  -- unknown -- functions $\mu_k$ which belong to the following space $\M$:
\be
\M:=\{\psi \in C^{0,\eta}([a,b]) \hbox{ s. t. }\psi \hbox{ is piecewise analytic on }(a,b)\},\label{hypmu}
\ee
for some $\eta \in (0,1]$. The space $C^{0,\eta}$ corresponds to H\"older continuous functions with exponent $\eta$ (see e.g. \cite{Fri64}).
A function $\psi \in C^{0,\eta}([a,b])$ is called piecewise
analytic if their exist $n>0$ and an increasing sequence $(\kappa_j)_{1\le j\le n}$ such that $\kappa_1=a$, $\kappa_n=b$, and $$\hbox{for all }x\in(a,b), \ \psi(x)=\sum_{j=1}^{n-1} \chi_{[\kappa_j,\kappa_{j+1})}(x)\varphi_j(x),$$for some analytic functions $\varphi_j$, defined on the intervals $[\kappa_j,\kappa_{j+1}]$, and where $\chi_{[\kappa_j,\kappa_{j+1})}$ are the characteristic functions of the intervals $[\kappa_j,\kappa_{j+1})$ for $j=1,\ldots, n-1$. In particular, if $\psi\in\M$, then, for each $x\in[a,b)$ (resp. $x\in(a,b]$), there exists $r=r_x>0$ such that $\psi$ is analytic on $[x,x+r]$ (resp. $[x-r,x]$).

The assumptions on the function $g$ are:
\be
g(\cdot,u) \in C^{0,\eta}([a,b]) \hbox{ for all }u\in \R,\ g(x,\cdot) \in C^{1}(\R) \hbox{ for all }x\in [a,b]\hbox{ and }g(\cdot,0)=0\hbox{ in }[a,b].\label{hyp:g}
\ee

We also assume that the diffusion coefficient $D$ is positive and that the boundary coefficients satisfy:
\be
 \alpha_1,\alpha_2,\beta_1,\beta_2\ge 0 \hbox{ with }\alpha_1+\beta_1>0 \hbox{ and }\alpha_2+\beta_2>0. \label{HypBC}
\ee

We furthermore make  the following hypotheses on the initial condition:
 \be u^0> 0 \hbox{ in }(a,b) \hbox{ and }u^0\in C^{2,\eta}([a,b]),\label{HypDI}\ee
that is $u^0$ is a $C^2$ function such that $(u^0) ''$ is H\"older continuous. In addition to that, we assume the following
compatibility conditions:
\be
\baa{rll}
     \alpha_1 u^0 (a) - \beta_1 (u^0) ' (a)=0\ \hbox{ and} & -D\,(u^0) ''(a)=g(a,0) & \!\!\hbox{if }\beta_1=0,\vspace{3pt}\\
     \alpha_2 u^0 (b) + \beta_2 (u^0) ' (b)=0\ \hbox{ and} & -D\,(u^0) ''(b)=g(b,0) & \!\!\hbox{if }\beta_2=0.
    \eaa
\label{HypDI2}
\ee

Under the assumptions (\ref{hypmu})-(\ref{HypDI2}), for each sequence $(\mu_k)_{1\leq k\leq N} \in \M^N$, there exists a time $T^{u^0}_{(\mu_k)}\in (0,+\infty]$ such that the problem $(\mathcal{P}^{u^0}_{(\mu_k)})$ has a unique solution
$\ds u \in C^{2,\eta}_{1,\eta/2}\left([0,T^{u^0}_{(\mu_k)})\times[a,b]\right)$ (i.e. the derivatives up to order two in $x$ and order one in $t$ are H\"older continuous). In the sequel, even if it means decreasing $T^{u^0}_{(\mu_k)}$ in some cases and dropping the indices $(\mu_k)$ and $u^0$, we only deal with values of $t$ smaller than $T$ so that the problem $(\mathcal{P}^{u^0}_{(\mu_k)})$ is well posed.
Existence, uniqueness and regularity of the solution $u$ are classical (see e.g.
\cite{Pao92}).

Our main result is a uniqueness result for the sequence of coefficients $(\mu_k)_{1\le k\le N}$  associated with observations of the solution and of its spatial derivative at a single point $x_0$ in $[a,b]$. Consider  $N$ ordered initial conditions $u^0_i$ and, for each sequence $(\mu_k)_{1\le k\le N}$, let $u_i$ be the solution of $(\mathcal{P}^{u^0_i}_{(\mu_k)}).$ Our result shows that for any $\varepsilon>0$ the function $$G : \baa{rll} \M^N & \to & C^1((0,\varepsilon))^{2\, N}\\(\mu_k)_{1\le k\le N} & \mapsto & (u_i(\cdot,x_0),\partial u_i/\partial x(\cdot,x_0))_{1\le i\le N}\eaa,$$
is one-to-one. In other words, we have the following theorem:

\begin{theorem}
\label{th:uniq2}
Let $N\in\N^*$, $(\mu_k)_{1\le k\le N}$ and $(\tmu_k)_{1\le k\le N}$ be two families of coefficients in $\M.$ Let~$\ds(u^0_i)_{1\le i\le N}$ be $N$ positive functions fulfilling (\ref{HypDI}) and (\ref{HypDI2}) and such that $u^0_i(x)\neq u^0_j(x)$ for all $ x\in (a,b)$ and
all $i\neq j.$ Let  $u_i$ and $\ub_i$ be the solutions of  $(\mathcal{P}^{u^0_i}_{(\mu_k)})$ and
$(\mathcal{P}^{u^0_i}_{(\tmu_k)}),$ respectively, on $[0,T)\times[a,b]$.
 We assume that
 $u_i$ and $\ub_i$ satisfy at some
$x_0\in [a,b],$ and for some $\varepsilon\in(0,T]$:
\be
\label{hyp:mainasummption}
\left\{
\baa{rll}
\ds u_i(t,x_0) & = & \ds \ub_i(t,x_0),\vspace{3pt}\\
\ds \frac{\partial u_i}{\partial x}(t,x_0) & = & \ds \frac{\partial \ub_i}{\partial x} (t,x_0),
\eaa \right.  \hbox{ for all }t\in(0, \varepsilon)\hbox{ and all }i\in\{1,\cdots,N\}.
\ee
 Then $\mu_k\equiv\tmu_k$ on $[a,b]$ for all $k\in\{1,\cdots,N\}$, and consequently $u_i\equiv \ub_i$ in $[0,T)\times[a,b]$ for all $i$.
\end{theorem}

The main result in \cite{RoqCri10} was a particular case of Theorem \ref{th:uniq2}. A similar conclusion was indeed proved in the case $N=1$ and for $g(x,u)=-\gamma \, u^2.$ In such case, the determination of one coefficient~$\mu_1(x)$ only requires the knowledge of the initial condition $u^0$ and of $(u(t,x_0),\partial u/\partial x(t,x_0))$ for $t\in (0,\varepsilon).$ When $N\ge 2,$ the above theorem requires more than the knowledge of the initial condition for the determination of the coefficients: we need a control on the initial condition, which enables to obtain $N$ measurements of the solution of $(\mathcal{P}^{u^0}_{(\mu_k)}),$ starting from $N$ different initial conditions. A natural question is whether the result of Theorem \ref{th:uniq2} remains true when the number of measurements is smaller than $N.$ In Section \ref{section_nonuniq}, we prove that the answer is negative in general.

\section{Proof of Theorem \ref{th:uniq2}}

For the sake of clarity, we begin with proving Theorem \ref{th:uniq2} in the particular case $N=2$ (the proof in the case $N=1$ would be similar to that of \cite{RoqCri10}, which was concerned with $g(x,u)=-\gamma\,u^2$). We then deal with the general case of problems~$(\mathcal{P}^{u^0_i}_{(\mu_k)})$ and
$(\mathcal{P}^{u^0_i}_{(\tmu_k)})$ with $N\ge 1$.

\subsection{Proof of Theorem \ref{th:uniq2}, case $N=2.$}\label{section_proof}

Let $(\mu_1,\tmu_1)$ and $(\mu_2,\tmu_2)$ be two pairs of coefficients in $\M.$ Let $u^0_1(x), u^0_2(x)$ be two functions verifying~(\ref{HypDI}) and (\ref{HypDI2}) and such that $u^0_1(x)\neq u^0_2(x)$ in $(a,b).$
Let $u_1$ and $\ub_1$ be respectively the solutions of $(\mathcal{P}^{u^0_1}_{\mu_1,\mu_2})$ and $(\mathcal{P}^{u^0_1}_{\tmu_1,\tmu_2})$ and $u_2$ and $\ub_2$ be the solutions of
$(\mathcal{P}^{u^0_2}_{\mu_1,\mu_2})$ and $(\mathcal{P}^{u^0_2}_{\tmu_1,\tmu_2}).$

We set, for $i=1,2,$
\[
   U_i:=u_i-\ub_i, \ m_1:=\mu_1-\tmu_1 \hbox{ and } m_2:=\mu_2-\tmu_2.
\]
The functions $U_i$ satisfy:
\begin{equation}
\frac{\partial U_i}{\partial t}- D  \frac{\partial^2 U_i}{\partial x^2}=b_i(t,x) U_i +h(x,u_i(t,x)),
\label{eqU}
\end{equation}
for $t \in [0,T)$ and $x\in[a,b],$ where
\be
\begin{array}{l}
 b_i(t,x)=\tmu_1(x)+\tmu_2(x)\left(u_i(t,x)+\ub_i(t,x)\right)+c_i(t,x),\vspace{3pt}\\
 c_i(t,x)= \left\{
           \begin{array}{l}
             \ds\frac{g(x,u_i(t,x))- g(x,\ub_i(t,x))}{u_i(t,x)-\ub_i(t,x)} \hbox{ if } u_i(t,x)\neq\ub_i(t,x),\vspace{3pt} \\
             \ds\frac{\partial g}{\partial u}(x,u_i(t,x))\hbox{ if } u_i(t,x)=\ub_i(t,x),
           \end{array}
          \right.\vspace{3pt}\\
 h(x,s)=s\left( m_1(x) + m_2(x)s\right),
\end{array}
\label{eq:b,h}
\ee
and the boundary and initial conditions:
\begin{equation}\left\{\baa{l}
\ds\alpha_1 U_i (t,a) - \beta_1 \frac{\partial U_i}{\partial x} (t,a)=0, \ t>0,\vspace{3pt} \\ [0.2 cm]
\ds\alpha_2 U_i (t,b) + \beta_2 \frac{\partial U_i}{\partial x} (t,b)=0, \ t>0,\vspace{3pt}\\
U_i(0,x)=0, \ x\in(a,b).
\eaa
\right.
\label{eqUb}
\end{equation}

\

Let us first assume that $x_0<b,$ and set:
\[
 \mathcal{A}_+=\Big\{x\geq x_0 \hbox{ s.t. } m_1(y)\equiv m_2(y)\equiv 0\hbox{ for all } y\in[x_0,x]\Big\},
\]
and
\[
x_1:=\left\{
\begin{array}{ll}
  \sup{(\mathcal{A}_+)} &\hbox{if}\ph \mathcal{A}_+\ph \hbox{is not empty},\vspace{3pt} \\
  x_0    &\hbox{if}\ph \mathcal{A}_+\ph \hbox{is empty}.
\end{array}
\right.
\]
If $x_1=b$, then $m_1(x)\equiv m_2(x) \equiv 0$ on $[x_0,b].$ Let us assume on the contrary that $x_1<b.$

\

\textit{Step 1: We show that there exist $\theta>0,$ $x_2 \in (x_1,b)$ and $j\in\{1,2\}$ such that the function
$(t,x)\mapsto h(x,u_j(t,x))$ has a constant strict sign on $[0,\theta]\times(x_1,x_2]$, i.e. $ h(x,u_j(t,x))>0$ or $h(x,u_j(t,x))<0$ in $[0,\theta]\times(x_1,x_2].$}\par
To do so, let us define, for all $x\in [x_1,b)$:
\begin{equation}
 z(x)=\left\{\begin{array}{ll}
              \ds -\frac{m_1(x)}{m_2(x)} & \hbox{if }m_2(x)\not=0,\vspace{3pt}\\
              \ds\lim_{y\to x^+}{-\frac{m_1(y)}{m_2(y)}} & \hbox{if } m_1(x)\!=\!m_2(x)\!=\!0 \hbox{ and } m_2(y)\!\not=\!0 \hbox{ in a right neighborhod of }x,\vspace{3pt}\\
               +\infty & \hbox{otherwise}.
             \end{array}
\right.\label{eq:z}
\end{equation}
Whenever $m_2(x)\not=0$, $z(x)$ is a root of the polynomial $h(x,\cdot).$ Notice also that the limit $\ds\lim_{
y\to x^+}{-\frac{m_1(y)}{m_2(y)}}$ in the second case of the definition of $z(x)$ is well defined since $m_1$ and $m_2$ are analytic on $[x,y]$ for $y-x>0$ small enough.

Since, $u_1^0(x_1)\neq u_2^0(x_1)$,  we have
$ |u_{j}^0(x_1)-z(x_1)|>0$ for some $j\in\{1,2\}.$
Moreover, there exists $\delta>0$ such that $x_1+\delta<b$ and
\be
 |u_{j}^0(x)-z(x)|\geq r>0 \hbox{ on }[x_1,x_1+\delta]\hbox{, for some }r>0.
\label{eq:dif_z}
\ee
Indeed, if $z(x_1)=\pm\infty$, we clearly have \eqref{eq:dif_z}. If $z(x_1)\neq \pm\infty$, from \eqref{eq:z}, $z$ is continuous in a right neighborhood of $x_1.$
As the function $u_{j}^0$ is also continuous in this neighborhood, we get \eqref{eq:dif_z}.

Moreover, $\ds u_{j}(t,x)\in C^{2,\eta}_{1,\eta/2}\left([0,T)\times[a,b]\right).$ This implies that $u_{j}$ is continuous at $(t,x)=(0,x_1).$ As a consequence, there
exists $\theta>0$ small enough so that $|u_{j}(t,x)-u_{j}^0(x_1)|\leq r/4$ in $[0,\theta]\times[x_1,x_1+\theta]$, whence
\be
|u_{j}(t,x)-u_{j}^0(x)|\leq \frac{r}{2} \hbox{ in } (t,x)\in[0,\theta]\times[x_1,x_1+\theta].
\label{eq:r2}
\ee

Finally, setting $\delta'=\min\{\theta,\delta\}$, we have, from (\ref{eq:dif_z}) and (\ref{eq:r2}):
\be
\ds|u_j(t,x)-z(x)| \geq  \left|\,|u_j(t,x)-u_j^0(x)|-|u_j^0(x)-z(x)|\,\right| \geq \frac{r}{2} > 0
\label{eq:r3}
\ee
for all $(t,x) \in[0,\theta]\times[x_1,x_1+\delta']$.

Now, the definition of $x_1$ and the piecewise analyticity of $m_1$ and $m_2$ imply that there exists $\delta''\in (0,\delta')$ such that:
\be
\hbox{for all } x\in(x_1,x_1+\delta''], \hbox{ the polynomial function }h(x,\cdot)\hbox{ verifies } h(x,\cdot)\not \equiv 0\hbox{ in }\R.
\label{eq:h}
\ee
Indeed, assume on the contrary that there is a decreasing sequence $y_n\to x_1$ such that $h(y_n,\cdot)\equiv0$ in $\R.$ Since the functions $h(y_n,\cdot)$
are polynomial, we get $m_1(y_n)=m_2(y_n)=0$ for all $n\in\N.$ Besides, as $m_1$ and $m_2$ belong to $\mathcal{M}$, for $n$ large enough, $m_1$ and $m_2$ are analytic on $[x_1,y_n]$,
which implies that $m_1\equiv m_2\equiv0$ on $[x_1,x_n]$ for $n$ large enough. This contradicts the definition of~$x_1$ and we get \eqref{eq:h}.

From the expression (\ref{eq:b,h}) of $h$ and using \eqref{eq:h}, we observe that for each $x\in(x_1,x_1+\delta'']$, either $s=0$ is the unique solution of $h(x,s)=0$ or $m_2(x)\not=0$ and the equation $h(x,s)=0$ admits exactly two solutions $s=0$ and $s=z(x).$
Let us set $x_2=x_1+\delta''\in (x_1,b).$ From the strong parabolic maximum principle $u_j(t,x)>0$ in $[0,\theta]\times(a,b).$ Thus, using \eqref{eq:r3} we finally get:
$$h(x,u_j(t,x))\not=0 \hbox{ in }[0,\theta]\times (x_1,x_2].$$
This concludes the proof of step $1.$

\

\textit{Step 2: We prove that $x_1=b$.}\par
From Step 1, let us assume in the sequel -- without loss of generality -- that $(t,x)\mapsto h(x,u_j(t,x))$ is positive on $[0,\theta]\times(x_1,x_2]$ (the case $h(x,u_j(t,x))<0$ could be treated similarly). Then, from the definition of $x_1,$ we deduce that $h(x,u_j(t,x))$ is nonnegative for all $(t,x)\in[0,\theta]\times[x_0,x_2].$

Since $h(x_2,u_j(0,x_2))>0$ and $U_j(0,\cdot)\equiv0$, it follows from (\ref{eqU}) that
\[
\frac{\partial U_j}{\partial t}(0,x_2)= h(x_2,u_j(0,x_2))>0.
 \]
Thus, for $\varepsilon'\in(0,\theta)$ small enough, $U_j(t,x_2)>0$ for  $t\in(0,\varepsilon')$. As a consequence, and from the assumption of Theorem \ref{th:uniq2}, $U_j$ satisfies:
\begin{equation*}
\left\{\baa{l}
\ds\frac{\partial U_j}{\partial t} - D  \frac{\partial^2 U_j}{\partial x^2} - b_j(t,x) U_j \ge 0, \ t \in (0, \varepsilon'), \ x\in [x_0,x_2],\vspace{3pt}
\\
U_j (t,x_0)=0 \hbox{ and }U_j(t,x_2)>0 , \ t\in(0,\varepsilon'),\vspace{3pt} \\
U_j(0,x)=0, \ x\in (x_0,x_2).
\eaa
\right.
\end{equation*}
Moreover, the weak and strong parabolic maximum principles give that
 $U_j(t,x)>0$ in $(0,\varepsilon')\times(x_0,x_2).$
Since $U_j(t,x_0)=0$, the Hopf's lemma  also implies that
$$\ds{\frac{\partial U_j}{\partial x}(t,x_0)>0}\hbox{ for all }t\in (0,\varepsilon').$$
This contradicts the assumption (\ref{hyp:mainasummption}) of Theorem  \ref{th:uniq2}. Finally, we necessarily have $x_1=b$ and therefore $m_1\equiv m_2\equiv0$ on $[x_0,b].$

\

\textit{Step 3: We prove that $m_1\equiv m_2\equiv 0$ on $[a,b]$.}\par
Assuming that $x_0>a$, and setting:
\[
 \mathcal{A}_-=\Big\{x\leq x_0 \hbox{ s.t. } m_1(y)=m_2(y)=0 \hbox{ for all }y\in [x,x_0]\Big\},
\]
and
\[
y_1:=\left\{
\begin{array}{ll}
  \inf{(\mathcal{A}_-)} &\hbox{if}\ph \mathcal{A}_-\ph \hbox{is not empty},\vspace{3pt}\\
  x_0    &\hbox{if}\ph \mathcal{A}_-\ph \hbox{is empty},
\end{array}
\right.
\]
we can prove, by applying the same arguments as above, that $y_1=a$ and consequently $m_1\equiv m_2\equiv0$ on $[a,x_0].$

Finally, $m_1\equiv m_2\equiv 0$ on $[a,b]$ which concludes the proof of Theorem \ref{th:uniq2} in the case $N=2$. \carre

\subsection{Proof of Theorem \ref{th:uniq2}, general case $N\ge 1.$}\label{section_proof2}

We set for all $i,k\in\{1,\cdots,N\},$
\[
U_i:=u_i-\ub_i, \ph m_k:=\mu_k-\tmu_k.
\]
The functions $U_i$ satisfy:
\begin{equation}
\frac{\partial U_i}{\partial t}- D  \frac{\partial^2 U_i}{\partial x^2}=b_i(t,x) U_i +h(x,u_i(t,x)),
\label{eqU2}
\end{equation}
for $t\in  [0,T)$ and $x\in[a,b],$ where
\bes
\begin{array}{l}
   \ds b_i(t,x)=\left\{
           \begin{array}{l}
             \ds\frac{(\tilde f+g)(x,u_i(t,x))-(\tilde f+g)(x,\ub_i(t,x))}{u_i(t,x)-\ub_i(t,x)} \hbox{ if } u_i(t,x)\neq\ub_i(t,x),\vspace{3pt} \\
             \ds\frac{\partial (\tilde f+g)}{\partial u}(x,u_i(t,x))\hbox{ if } u_i(t,x)=\ub_i(t,x),
           \end{array}
          \right.\vspace{3pt}\\
 \ds h(x,s)=\sum_{k=1}^{N}{m_k(x)s^k},
\end{array}
\ees
and
$$\ds \tilde f(x,u)=\sum_{k=1}^{N}{\tmu_k(x)u^k}.$$
Moreover, the functions $U_i$ satisfy the following boundary and initial conditions:
\begin{equation}\left\{\baa{l}
\ds\alpha_1 U_i (t,a) - \beta_1 \frac{\partial U_i}{\partial x} (t,a)=0, \ t>0,\vspace{3pt}\\ [0.2cm]
\ds\alpha_2 U_i (t,b) + \beta_2 \frac{\partial U_i}{\partial x} (t,b)=0, \ t>0,\vspace{3pt}\\
U_i(0,x)=0, \ x\in(a,b).
\eaa
\right.
\label{eqUb2}
\end{equation}

\

Let us set
\[
 \mathcal{A}_+=\Big\{x\geq x_0 \hbox{ s.t. } m_k\equiv 0 \hbox{ on } [x_0,x]\hbox{ for all }k\in\{1,\cdots,N\}\Big\},
\]
and
$$
x_1:=\left\{
\begin{array}{ll}
  \sup{(\mathcal{A}_+)} &\hbox{if}\ph \mathcal{A}_+\ph \hbox{is not empty}, \\
  x_0    &\hbox{if}\ph \mathcal{A}_+\ph \hbox{is empty}.
\end{array}
\right.
$$
If $x_1=b$, then for all $k$, $m_k(x)\equiv 0$ on $[x_0,b].$ Let us assume by contradiction that $x_1<b.$ As in the case $N=2,$ we prove, as a first step, that
there exist $\theta\in(0,T),$ $x_2 \in (x_1,b)$ and $j\in\{1,\ldots,N\}$ such that the function
$(t,x)\mapsto h(x,u_j(t,x))$ has a constant strict sign on $[0,\theta]\times(x_1,x_2]$, i.e. $ h(x,u_j(t,x))>0$ or $h(x,u_j(t,x))<0$ in $[0,\theta]\times(x_1,x_2].$

To do so, observe that, from the definitions of $x_1$ and $\M$, there exists $\delta>0$ such that $x_1+\delta<b$ and all functions $m_k$ are analytic on $[x_1,x_1+\delta]$ and not all identically zero. Therefore, the integer
$$\rho=\max\Big\{\rho'\in\N,\ m_k(x)=O\Big((x-x_1)^{\rho'}\Big)\hbox{ as }x\to x_1^+\hbox{ for all }k\in\{1,\ldots,N\}\Big\}$$
is well-defined. Furthermore, the function $h$ can then be written as
$$h(x,s)=(x-x_1)^{\rho}\times H(x,s)\hbox{ for all }(x,s)\in[x_1,x_1+\delta]\times\R,$$
where
$$H(x,s)=M_1(x)\,s+\cdots+M_N(x)\,s^N$$
and the functions $M_1,\ldots,M_N$ are analytic on $[x_1,x_1+\delta]$ and not all zero at the point $x_1$ (namely, there exists $i\in\{1,\ldots,N\}$ such that $M_i(x_1)\neq 0$). In other words, the polynomial $H(x_1,\cdot)$ is not identically zero. Since its degree is not larger than $N$ and since $H(x_1,0)=0$ and the real numbers $u^0_1(x_1),\ldots,u^0_N(x_1)$ are all positive and pairwise different, there exists $j\in\{1,\ldots,N\}$ such that
$$H(x_1,u^0_j(x_1))\neq 0.$$
By continuity of $H$ in $[x_1,x_1+\delta]\times\R$ and of $u_j$ on $[0,T)\times[a,b]$, it follows that there exist $\theta\in(0,T)$ and $x_2\in(x_1,b)$ such that
$$H(x,u_j(t,x))\neq 0\hbox{ for all }(t,x)\in[0,\theta]\times[x_1,x_2].$$
Consequently,
$$h(x,u_j(t,x))\neq 0\hbox{ for all }(t,x)\in[0,\theta]\times(x_1,x_2].$$

The remaining part of the proof of Theorem \ref{th:uniq2} in the general case $N\ge 1$ is then similar to Steps 2 and 3 of the proof in the particular case $N=2.$ Namely, we eventually get a contradiction with the assumption that $\frac{\partial U_j}{\partial x}(t,x_0)=0$ for all $t\in(0,\epsilon)$, yielding $x_1=b$ and $m_k=0$ on $[x_0,b]$ for all $k\in\{1,\ldots,N\}$. Similarly, $m_k=0$ on $[a,x_0]$ for all $k\in\{1,\ldots,N\}$. $\Box$

\section{Non-uniqueness results}\label{section_nonuniq}

This section deals with non-uniqueness results for the coefficients $(\mu_k)$ in $(\mathcal{P}^{u^0}_{(\mu_k)})$ under assumptions weaker than those of Theorem \ref{th:uniq2}. These results emphasize the optimality of the assumptions of Theorem \ref{th:uniq2}.

\

\textit{1--Number of measurements smaller than number of unknown coefficients.}

\

We give a counter-example to the uniqueness result of Theorem \ref{th:uniq2} in the case where the number of measurements is smaller than the number of unknown coefficients $N$.

Assume that the coefficients $\mu_1,\ldots, \mu_N$ are constant, not all zero, and such that the polynomial $$\ds{f(x,u)=f(u)=\sum_{k=1}^N\, \mu_k  \, u^k}$$admits  exactly $N-1$ positive and distinct roots $z_1,\ldots,z_{N-1}$. Assume furthermore that $\alpha_1=\alpha_2=0$ (Neumann boundary conditions) and that $g\equiv 0$. Then for each $i=1,\ldots,N-1$, $z_i$ is a (stationary) solution of $(\mathcal{P}^{z_i}_{(\mu_k)}).$ Consider a similar problem with the coefficients $\tilde \mu_k= \tau \, \mu_k$ for $\tau \neq 1$ and $k=1,\ldots,N.$ Then, again, for each
$i=1,\ldots,N-1$, $z_i$ is a solution of $(\mathcal{P}^{z_i}_{(\tmu_k)}).$ In particular, assumption \eqref{hyp:mainasummption} is fulfilled at any point $x_0\in [a,b]$ for $k\in\{1,\ldots,N-1\}.$ However $(\mu_1(x),\ldots,\mu_N(x))\neq (\tilde \mu_1(x),\ldots,\tilde \mu_N(x)).$

This shows that the determination of $N$ coefficients $(\mu_k)_{1\le k\le N}$ requires in general $N$ observations of the solution of $(\mathcal{P}^{u^0}_{(\mu_k)}),$ starting from $N$ different initial conditions.

\

\textit{2--Lack of measurement of the spatial derivatives.}

\

We show that if hypothesis \eqref{hyp:mainasummption} in Theorem \ref{th:uniq2} is replaced with the weaker assumption:
\be
\label{hyp:mainasummptionbis}
u_i(t,x_0)  =  \ub_i(t,x_0), \hbox{ for all }t\in(0, \varepsilon)\hbox{ and all }i\in\{1,\cdots,N\},
\ee
then the conclusion of the theorem is false in general.

Let $(\mu_k)_{1\le k\le N} \in \M^N$ and assume that $\alpha_1=\alpha_2=0$ (Neumann boundary conditions). Let $\ds(u^0_i)_{1\le i\le N}$ satisfy the assumptions of Theorem \ref{th:uniq2} and assume furthermore that  the functions $u^0_i$ and $g(\cdot,u)$ are symmetric with respect to $x=(a+b)/2,$ i.e.
$$\left\{\begin{array}{l}
\ds(u^0_i)_{1\le i\le N}(x)=\ds(u^0_i)_{1\le i\le N}(b-(x-a))\vspace{3pt}\\
g(x,\cdot)=g(b-(x-a),\cdot)\end{array}\right.\hbox{ for all }x\in [a,b].$$
Let $\tmu_k:=\mu_k(b-(x-a))$ for all $x\in [a,b]$ and $k\in\{1,\cdots,N\}.$

Then, we claim that the solutions
$u_i$ and $\ub_i$ of $(\mathcal{P}^{u_i^0}_{(\mu_k)})$ and $(\mathcal{P}^{u_i^0}_{(\tmu_k)})$ satisfy \eqref{hyp:mainasummptionbis} at $x_0=\frac{a+b}{2}$ and for $\varepsilon$ small enough. Indeed, we observe that, for each $i\in\{1,\cdots,N\},$ $\ub_i(t,b-(x-a))$ is a solution of $(\mathcal{P}^{u_i^0}_{(\mu_k)}).$ By uniqueness, we have $$u_i(t,x)=\ub_i(t,b-(x-a)), \hbox{ for all }t\in (0,T),  \hbox{ all }x\in [a,b]  \hbox{ and all }i\in\{1,\cdots,N\}.$$In particular,  $u_i(t,\frac{a+b}{2})  =  \ub_i(t,\frac{a+b}{2})$ for $t\in (0,T)$ and $i\in\{1,\cdots,N\}.$

This shows that the assumption \eqref{hyp:mainasummptionbis} alone is not sufficient to determine the coefficients $(\mu_k)_{1\le k\le N}.$

The above result is an adaptation of Proposition 2.3 in \cite{RoqCri10} to the general case $N\ge 1.$

\

\textit{3--Time-dependent coefficients.}

\

We show here that the result of Theorem \ref{th:uniq2} is not true in general when the coefficients $(\mu_k)$ are allowed to depend on the variable $t.$

We place ourselves in the simple case $N=1$ and $g\equiv 0$, and we assume that $\alpha_1=\alpha_2=0$ (Neumann boundary conditions). We assume that $(a,b)=(0,\pi),$ and we set  $u(t,x)=1+t\, \cos^2(x)$ and $\ub(t,x)=1+t\, \sin^2(2x).$

Let us set $$\mu_1(t,x)=\frac{1}{u}\left(\frac{\partial u}{\partial t}-D\, \frac{\partial^2 u}{\partial x^2}\right)(t,x) \hbox{ and }\tmu_1(t,x)=\frac{1}{\ub}\left(\frac{\partial \ub}{\partial t}-D\, \frac{\partial^2 \ub}{\partial x^2}\right)(t,x)),$$i.e.,
$$\mu_1(t,x)=\frac{(4Dt+1)\cos^2(x)-2Dt}{1+t\, \cos^2(x)} \hbox{ and }\tmu_1(t,x)=\frac{(16Dt+1)\sin^2(2x)-8Dt}{1+t\, \sin^2(2x)}.$$
Then, for each $t\in [0,+\infty),$ $\mu_1(t,\cdot)$ and $\tmu_1(t,\cdot)$ belong to $\M.$ The functions $u$ and $\ub$ are solutions of $(\mathcal{P}^{1}_{(\mu_1)})$ and $(\mathcal{P}^{1}_{(\tmu_1)}),$ respectively, and they satisfy the assumption \eqref{hyp:mainasummption} of Theorem \ref{th:uniq2} at $x_0=\pi/2.$ However, the conclusion of Theorem \ref{th:uniq2} does not hold since $\mu_1 \not \equiv \tmu_1.$

\

\textit{4--Unknown initial data.}

\

We show here that the result of Theorem \ref{th:uniq2} is not true in general if the functions $u_i$ and $\ub_i$ are solutions of $(\mathcal{P}^{u^0_i}_{(\mu_k)})$ and $(\mathcal{P}^{\ub^0_i}_{(\tmu_k)}),$  with $u^0_i\not \equiv \ub^0_i.$ This means that the coefficients $(\mu_k)$ cannot be determined, given only the measurements $(u_i(t,x_0),\partial u_i/\partial x(t,x_0))_{1\le i\le N}$ for $t\in (0,\varepsilon),$ if the initial conditions $u^0_i$ are unknown.

We build an explicit counter-example in the  simple case $N=1$ and $g\equiv 0,$ with $\alpha_1=\alpha_2=0$. Assume that $(a,b)=(0,\pi),$ and
let us set $u(t,x)=(1+\cos^2(x))\, e^{\rho \, t}$ for some $\rho>0,$ and $\ub(t,x)=(1+\sin^2(2x))\, e^{\rho \, t}.$ We furthermore set $$\mu_1(x)=\frac{1}{u}\left(\frac{\partial u}{\partial t}-D\, \frac{\partial^2 u}{\partial x^2}\right) \hbox{ and }\tmu_1(x)=\frac{1}{\ub}\left(\frac{\partial \ub}{\partial t}-D\, \frac{\partial^2 \ub}{\partial x^2}\right),$$i.e.,
$$\mu_1(x)=\frac{(4D+\rho)\cos^2(x)+\rho-2\, D}{1+\cos^2(x)} \hbox{ and }\tmu_1(x)=\frac{(16D+\rho)\sin^2(2x)+\rho-8\, D}{1+\sin^2(2x)},$$thus $\mu_1$ and $\tmu_1$ belong to $\M.$ Besides, $u$ and $\ub$ are solutions of
$(\mathcal{P}^{1+\cos^2(x)}_{\mu_1})$ and $(\mathcal{P}^{1+\sin^2(2x)}_{\tmu_1})$ respectively, and they satisfy the assumption \eqref{hyp:mainasummption} of Theorem \ref{th:uniq2} at $x_0=\pi/2.$ However, we obviously have $\mu_1\not \equiv \tmu_1.$

\section{Numerical determination of several coefficients}\label{section_num}

In the particular case $N=1$, it was shown in \cite{RoqCri10} that the measurements
 \eqref{hyp:mainasummption} of Theorem \ref{th:uniq2} are sufficient to obtain a good numerical approximation of a coefficient $\mu_1(x).$
In this section, we check whether the measurements \eqref{hyp:mainasummption} of Theorem \ref{th:uniq2} also allow for an accurate reconstruction of $N$ coefficients $(\mu_k)_{1\le k \le N}$ in the cases $N=2$ and $N=3.$

Given the initial data $u^0_i$ and the measurements $u_{i}(t ,x_0)$ and $\frac{\partial u_i}{\partial x}(t,x_0),$ for $t\in (0,\varepsilon)$ and $i\in\{1,\cdots,N\},$  we can look for  the sequence  $(\mu_k)_{1\le k \le N}$ as a minimizer of some functional $G_{(\mu_k)}.$ Indeed, for any
 sequence $(\tmu_k)_{1\le k \le N}$ in $\M^N,$ the distance between the measurements of the solutions $u_i$ of $(\mathcal{P}^{u^0_i}_{(\mu_k)})$ and $\ub_i$ of $(\mathcal{P}^{u^0_i}_{(\tmu_k)}),$ $i\in\{1,\cdots,N\}$ can be evaluated through the function:
$$G_{(\mu_k)}[(\tmu_k)]=\sum_{i=1}^N \|u_{i}(\cdot ,x_0)-\ub_{i}(\cdot,x_0)\|_{L^2(0,\varepsilon)}+\|\frac{\partial u_i}{\partial x}(\cdot,x_0)-\frac{\partial \ub_i}{\partial x}(\cdot,x_0)\|_{L^2(0,\varepsilon)}.$$ Then, $G_{(\mu_k)}[(\mu_k)]=0$ and from Theorem \ref{th:uniq2} this is the unique global minimum of $G_{(\mu_k)}$ in $\M^N$.

In our numerical computations,  we fixed  $(a,b)=(0,1),$ $D=0.1$, $\alpha_1=\alpha_2=0$ and $\beta_1=\beta_2=1$ (Neumann boundary conditions),  $x_0=2/3$ and $\varepsilon=0.3.$ Besides, we assumed the coefficients $\mu_k$ to  belong to  a finite-dimensional subspace $E\subset \M:$
\begin{equation*} \label{eq:defE} E:=\left\{\rho \in  C^{0,\eta}([0,1]) \, | \, \exists \  (h_j)_{0\le j \le n}\in \R^{n+1}, \ \rho(x)=\sum_{j=0}^{n}h_j \cdot J\left((n-2)\left(x-c_j\right)\right) \hbox{ on }[0,1]\right\},\end{equation*} with $c_j=\frac{j-1}{n-2}$ and
$J(x)=\ds \frac{(x-2)^4(x+2)^4}{2^8}$ if $x\in(-2,2),$  and $J(x)= 0$  otherwise. In our computations, the integer $n$ was set to $10.$

\begin{figure}
\centering
\subfigure[Determination of 2 coefficients]{
\includegraphics*[width=8cm]{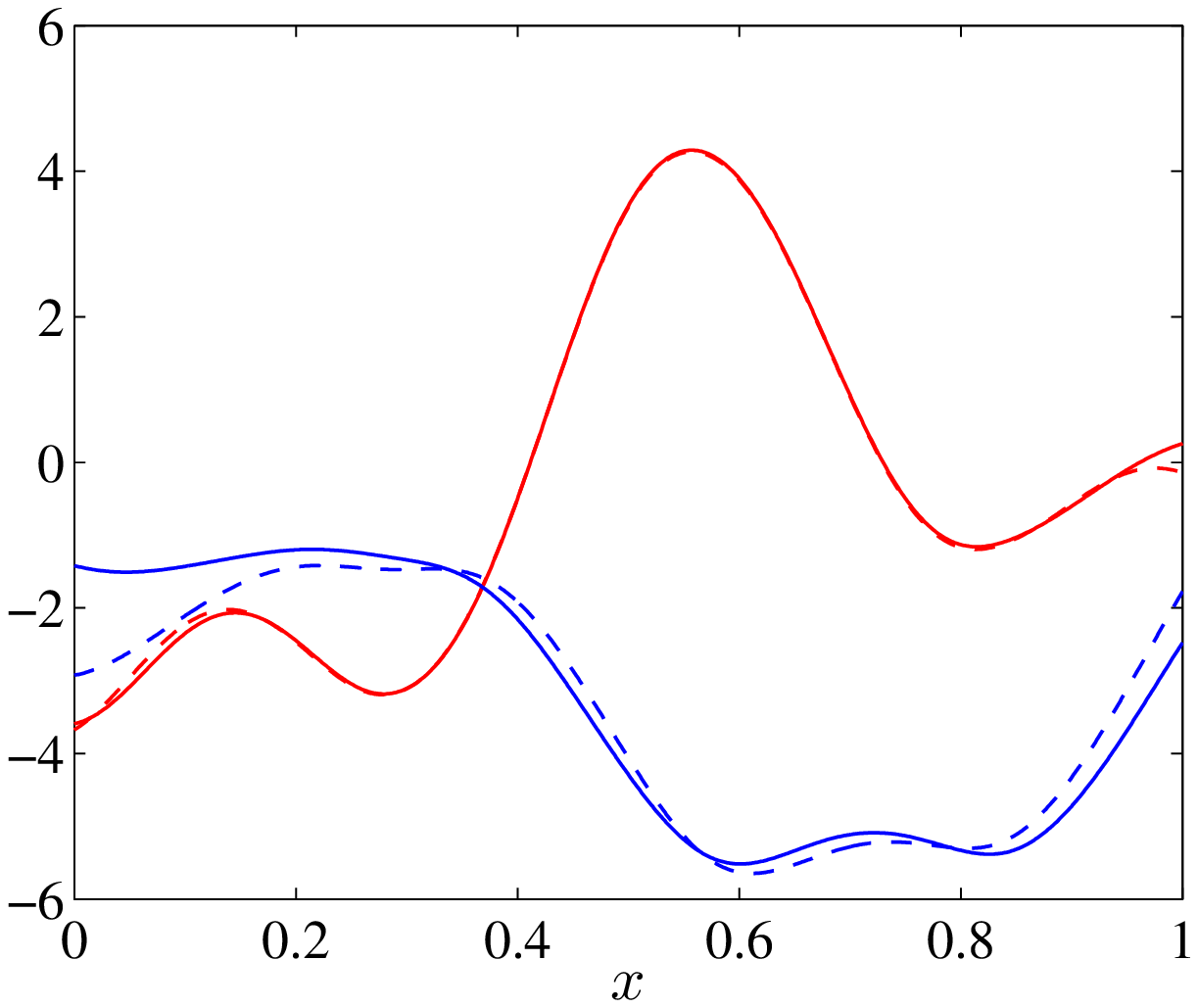}}
\subfigure[Determination of 3 coefficients]{\includegraphics*[width=8cm]{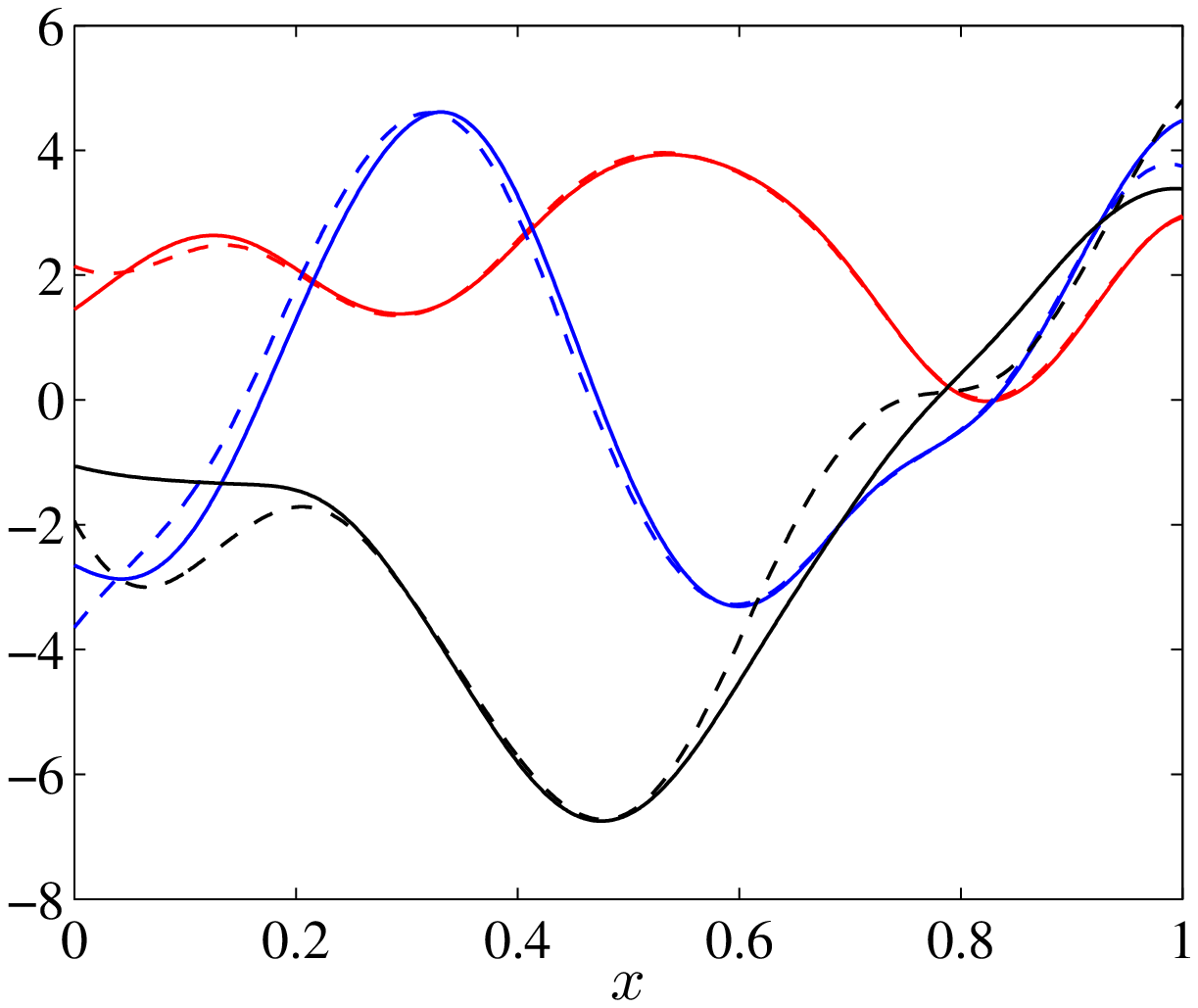}}
\caption{ (a) Plain  lines: examples of functions $\mu_1$ (red line) and $\mu_2$ (blue line) in $E.$ Dashed lines:
 the functions $\mu^*_1$ (in red) and $\mu^*_2$ (in blue)  obtained by minimizing $G_{(\mu_1,\mu_2)}.$ In this case $\ds{\|\mu_1-\mu^*_1\|_{L^2(0,1)}+\|\mu_2-\mu^*_2\|_{L^2(0,1)}=0.15},$ and $G_{(\mu_1,\mu_2)}[(\mu_1^*,\mu_2^*)]=9 \cdot 10^{-6}.$ (b) Plain  lines: functions $\mu_1$ (red line), $\mu_2$ (blue line) and $\mu_3$ (black line) in $E.$ Dashed lines:
 the functions $\mu^*_1$ (in red), $\mu^*_2$ (in blue), $\mu^*_3$ (in black) obtained by minimizing $G_{(\mu_1,\mu_2,\mu_2)}.$ Here, $\ds{\|\mu_1-\mu^*_1\|_{L^2(0,1)}+\|\mu_2-\mu^*_2\|_{L^2(0,1)}+\|\mu_3-\mu^*_3\|_{L^2(0,1)}}=0.38$ and $G_{(\mu_1,\mu_2,\mu_3)}[(\mu_1^*,\mu_2^*,\mu_3^*)]=3 \cdot 10^{-5}.$
 }
\label{fig:exple}
\end{figure}

\underline{Case $N=2:$} $25$ couples of functions $(\mu_1,\mu_2)$  have been randomly sampled in $E^2:$ for $k=1$ and $k=2,$
the components $h_j^k$, in the expression $$\mu_k(x)=\sum_{j=0}^{n}h_j^k \cdot J\left((n-2)\left(x-c_j\right)\right),$$were randomly drawn from  a uniform distribution in $(-5,5).$

Starting from the initial data $u_1^0\equiv 0.1,$ and $u_2^0\equiv 0.2,$ the corresponding values of $u_{1}(t ,x_0),$  $\frac{\partial u_1}{\partial x}(t,x_0),$ $u_{2}(t ,x_0),$  $\frac{\partial u_2}{\partial x}(t,x_0)$ were recorded\footnote{Numerical computation of $u$ and $\ub$ were carried out with Comsol Multiphysics$^{\circledR}$ time-dependent solver. We used a second order finite element method (FEM) with 960 elements. This solver uses a method of lines approach incorporating variable order and variable stepsize backward differentiation formulas.}, which enabled us to compute $G_{(\mu_1,\mu_2)}[(\tmu_1,\tmu_2)]$ for any couple $(\tmu_1,\tmu_2)$ in  $E^2$. The minimizations\footnote{The minimizations of the functions $G_{(\mu_k)}$ were performed using MATLAB's$^{\circledR}$ \textit{fminunc} solver. This optimization algorithm uses a Quasi-Newton method with a mixed quadratic and cubic line search procedure. The stopping criterion was based on the number of evaluations of the function $G_{(\mu_1,\mu_2)}$ which was set at $4\cdot 10^3.$} of the functions $G_{(\mu_1,\mu_2)}$ lead to
$25$ couples  $(\mu_1^*,\mu_2^*)$, each one corresponding to a computed  minimizer of the function $G_{(\mu_1,\mu_2)}$.

The average value of the quantity $\ds{\|\mu_1-\mu^*_1\|_{L^2(0,1)}+\|\mu_2-\mu^*_2\|_{L^2(0,1)}},$ over the $25$ samples of couples $(\mu_1,\mu_2)$  is $0.25$. The corresponding average value of $G_{(\mu_1,\mu_2)}(\mu_1^*,\mu_2^*)$ is $1.5 \cdot 10^{-5}.$
Fig. \ref{fig:exple} (a) depicts an example of a couple $(\mu_1,\mu_2)$ in $E$, together with the couple $(\mu_1^*,\mu_2^*)$ which was obtained by minimizing $G_{(\mu_1,\mu_2)}.$

\underline{Case $N=3:$} in this case, the minimization of the function $G_{(\mu_1,\mu_2,\mu_3)}$ is more time-consuming. We therefore focused on a unique example of a triple $(\mu_1,\mu_2,\mu_3)$ in $E^3.$ The initial data were chosen as follows: $u_1^0\equiv 0.1,$ $u_2^0\equiv 0.2,$ and $u_3^0\equiv 0.3.$ Fig. \ref{fig:exple} (b) depicts the triple $(\mu_1,\mu_2,\mu_3)$ in $E$, together with the triple $(\mu_1^*,\mu_2^*,\mu_3^*)$ obtained by minimizing $G_{(\mu_1,\mu_2,\mu_3)}.$

\section{Discussion}

We have obtained a uniqueness result in the inverse problem of determining several non-constant coefficients of reaction-diffusion equations. With a reaction term containing an unknown polynomial part of the form $\sum_{k=1}^N\, \mu_k(x) \, u^k,$ our  result provides a sufficient condition for the uniqueness of the determination of this nonlinear polynomial part.

This sufficient condition, which is detailed in Theorem \ref{th:uniq2}, involves pointwise measurements of the solution $u(t,x_0)$  and of its spatial derivative $\partial u / \partial x(t,x_0)$ at a single point $x_0,$ during a time interval $(0,\varepsilon),$ and starting with $N$ nonintersecting initial conditions.

The results of Section \ref{section_nonuniq} show that most conditions of  Theorem \ref{th:uniq2} are in fact necessary. In particular,
the first counter-example of Section \ref{section_nonuniq} shows that, for the result of Theorem \ref{th:uniq2} to hold in general, the number of measurements of the couple $(u,\partial u / \partial x)(t,x_0)$ needs to be at least equal to the degree ($N$) of the unknown polynomial term.

From a practical point of view, such measurements can be obtained if one has a control on the initial condition. Nevertheless, since our result does not provide a stability inequality, the possibility to do a numerical reconstruction of the unknown coefficients $\mu_k,$ on the basis of pointwise measurements,  was uncertain. In Section \ref{section_num}, we have shown in the cases $N=2$ and $N=3$ -- which include the classical models \eqref{eq_IntroSke} and \eqref{bistable_heter} -- that such measurements can indeed lead to good numerical approximations of the unknown coefficients, at least if they are assumed to belong to a known finite-dimensional space.

\section*{Acknowledgements}

The authors are supported by the French ``Agence Nationale de la Recherche"
within the projects ColonSGS (third and fourth authors), PREFERED (second, third and fourth authors) and URTICLIM (second and fourth authors). The third author is also indebted to the Alexander von~Humboldt Foundation for its support.

\providecommand{\bysame}{\leavevmode\hbox to3em{\hrulefill}\thinspace}
\providecommand{\MR}{\relax\ifhmode\unskip\space\fi MR }
\providecommand{\MRhref}[2]{%
  \href{http://www.ams.org/mathscinet-getitem?mr=#1}{#2}
}
\providecommand{\href}[2]{#2}

\end{document}